\documentclass{amsart}
\usepackage{amsmath}
\usepackage{paralist}
\usepackage{amsfonts}
\usepackage{amssymb}
\usepackage{amsthm}
\usepackage{amscd}
\usepackage{float}
\usepackage{tikz}
\usepackage{graphicx}
\usepackage[colorlinks=true]{hyperref}
\hypersetup{urlcolor=blue, citecolor=red}
\usepackage{hyperref}

\newtheorem{theorem}{Theorem}[section]
\newtheorem{corollary}[theorem]{Corollary}

\newtheorem{lemma}[theorem]{Lemma}
\newtheorem{proposition}[theorem]{Proposition}

\theoremstyle{definition}

\newtheorem{remark}[theorem]{Remark}

\makeatletter
\@namedef{subjclassname@2020}{%
  \textup{2020} Mathematics Subject Classification}
\makeatother

\begin{document}
\title[A priori estimates for dNLS]{A priori estimates for the derivative nonlinear Schr\"{o}dinger equation}
\author{Friedrich Klaus}
\author{Robert Schippa}
\email{friedrich.klaus@kit.edu}
\email{robert.schippa@kit.edu}
\address{Fakult\"{a}t f\"{u}r Mathematik, Karlsruher Institut f\"{u}r Technologie, Englerstrasse 2, 76131 Karlsruhe, Germany}

\begin{abstract}
We prove low regularity a priori estimates for the derivative nonlinear Schr\"{o}dinger equation in Besov spaces with positive regularity index conditional upon small $L^2$-norm. This covers the full subcritical range. We use the power series expansion of the perturbation determinant introduced by Killip--Vi\c{s}an--Zhang for completely integrable PDE. This makes it possible to derive low regularity conservation laws from the perturbation determinant.
\end{abstract}

\maketitle
\section{Introduction}
In this note the following derivative nonlinear Schr\"odinger equation (\emph{dNLS}) is considered
\begin{equation}
\label{eq:dNLS}
\left\{ \begin{array}{cl}
i \partial_t q + \partial_{xx} q + i \partial_x(|q|^2 q)  &= 0 \quad (t,x) \in \mathbb{R} \times \mathbb{K}, \\
q(0) &= q_0 \in H^s(\mathbb{K}),
\end{array} \right.
\end{equation}
where $\mathbb{K} \in \{ \mathbb{R}, \mathbb{T} = (\mathbb{R} / (2 \pi \mathbb{Z})) \}$. In the seventies \eqref{eq:dNLS} was proposed as a model in plasma physics in \cite{Rogister1971,MioOginoMinamiTakeda1976,Mjolhus1976}.

In the following let $\mathcal{S}(\mathbb{R})$ denote the Schwartz functions on the line and $\mathcal{S}(\mathbb{T})$ smooth functions on the circle.
Here we prove a priori estimates
\begin{equation*}
\sup_{t \in \mathbb{R}} \Vert q(t) \Vert_{H^s} \lesssim_{s} \Vert q_0 \Vert_{H^s}, \quad 0 < s < \frac{1}{2},
\end{equation*}
where $q \in C^\infty(\mathbb{R}; \mathcal{S}(\mathbb{K}))$ is a smooth global solution to \eqref{eq:dNLS}, which is also rapidly decaying in the line case, \emph{conditional upon small} $L^2$\emph{-norm}. These estimates are the key to extend local solutions globally in time. Local well-posedness, i.e., existence, uniqueness and continuous dependence locally in time, in $H^{1/2}$ was proved by Takaoka \cite{Takaoka1999} on the real line and Herr \cite{Herr2006} on the circle. They proved local well-posedness via the contraction mapping principle, that is \emph{perturbatively}. Furthermore, they showed that the data-to-solution mapping fails to be $C^3$ below $H^{1/2}$ in either geometry, respectively. Moreover, Biagioni--Linares \cite{BiagioniLinares2001} showed that the data-to-solution mapping even fails to be locally uniformly continuous on the real line below $H^{1/2}$. Thus, the results on local well-posedness in $H^{1/2}$ are the limit of proving local well-posedness via fixed point arguments. However, on the real line \eqref{eq:dNLS} admits the scaling symmetry
\begin{equation}
\label{eq:Scaling}
q(t,x) \rightarrow \lambda^{-1/2} q(\lambda^{-2} t, \lambda^{-1} x),
\end{equation}
which distinguishes $L^2$ as scaling critical space. Hence, we still expect a milder form of local well-posedness in $H^s$ for $0 \leq s < 1/2$. By short-time Fourier restriction, Guo \cite{Guo2011} proved a priori estimates for $s>1/4$ on the real line, which the second author \cite{Schippa2017} extended to periodic boundary conditions.\\
Moreover, Gr\"{u}nrock \cite{Gruenrock2005} showed local well-posedness on the real line in Fourier Lebesgue spaces, which scale like $H^s$, $s>0$. Deng \emph{et al.} \cite{DengNahmodYue2019} recently extended this to periodic boundary conditions; see also the previous work \cite{GruenrockHerr2008}. 

Less is known about global well-posedness. Conserved quantities of the flow include the mass, i.e., the $L^2$-norm,
\begin{equation*}
M[q] = \int_{\mathbb{K}} |q|^2 dx,
\end{equation*}
the momentum, related with the $H^{1/2}$-norm,
\begin{equation*}
P[q] = \int_{\mathbb{K}} \mathrm{Im}(\bar{q} q_x) - \frac{1}{2} |q|^4 dx,
\end{equation*}
and the energy, related with the $H^1$-norm,
\begin{equation*}
E[q] = \int_{\mathbb{K}} |q_x|^2 - \frac{3}{2} |q|^2 \mathrm{Im}( \bar{q} q_x) + \frac{1}{2} |q|^6 dx.
\end{equation*}
A local well-posedness result in $L^2$ seems to be very difficult due to the scaling criticality. On the other hand, it is not straight-forward to use the other quantities to prove a global result due to lack of definiteness. The remedy in previous works was to impose a smallness condition on the $L^2$-norm and use the sharp Gagliardo-Nirenberg inequality.

Wu \cite{Wu2015} observed in the line case that combining several conserved quantities improves the $L^2$-threshold, which can be derived from the energy (cf. \cite{Wu2013}). Mosincat--Oh carried out the corresponding argument on the torus \cite{MosincatOh2015}. Additionally making use of the $I$-method (cf. \cite{CollianderKeelStaffilaniTakaokaTao2002, MiaoWuXu2011}), Guo--Wu \cite{GuoWu2017} proved global well-posedness in $H^{1/2}(\mathbb{R})$ for $\Vert u_0 \Vert^2_{L^2} < 4 \pi$, and Mosincat \cite{Mosincat2017} proved global well-posedness in $H^{1/2}(\mathbb{T})$ under the same $L^2$-smallness condition. The question of global well-posedness for arbitrary $L^2$-norm was still open at the time of the first submission of this paper. Afterwards, there were several new contributions to the global well-posedness of dNLS (\cite{TangXu2020,BahouriPerelman2020,KillipNtekoumeVisan2021,IsomMantzavinosStefanov2020}). Among these, Bahouri--Perelman \cite{BahouriPerelman2020} showed global well-posedness in $H^{1/2}(\mathbb{R})$ without smallness assumption on the $L^2$-norm. The new works are discussed at the end of the Introduction. Previously, Nahmod \emph{et al.} \cite{NahmodOhReyBelletStaffilani2012} proved a probabilistic global well-posedness result in Fourier Lebesgue spaces scaling like $H^{1/2-\varepsilon}(\mathbb{T})$. On the half-line and intervals endowed with Dirichlet boundary conditions, Wu \cite{Wu2013} and Tan \cite{Tan2004} showed the existence of finite-time blow-up solutions.

Kaup--Newell \cite{KaupNewell1978} already observed shortly after the proposal of \eqref{eq:dNLS} that it admits a Lax pair with operator
\begin{equation}
\label{eq:LaxOperator}
L(t;q) = 
\begin{pmatrix}
\partial + i \kappa^2 & - \kappa q \\
 \kappa \bar{q} & \partial - i \kappa^2
\end{pmatrix}
.
\end{equation}
Consequently, there are infinitely many conserved quantities of the flow. However, to the best of the authors' knowledge, there were no prior works using the complete integrability for solutions in \emph{unweighted} $L^2$-based Sobolev space, i.e., without imposing additional spatial decay. In particular, there were no results for periodic boundary conditions making use of the complete integrability before the present ones.

Via inverse scattering, Lee \cite{Lee1983,Lee1989} proved global existence and uniqueness for certain initial data in $\mathcal{S}(\mathbb{R})$. Later, Liu \cite{Liu2017} considered \eqref{eq:dNLS} with initial data in weighted Sobolev spaces $H^{2,2}(\mathbb{R})$ and proved global well-posedness via inverse scattering. See the subsequent works \cite{JenkinsLiuPerrySulem2018I,JenkinsLiuPerrySulem2018II} due to Jenkins \emph{et al.} for results addressing soliton resolution in weighted Sobolev spaces and \cite{JenkinsLiuPerrySulem2020Survey} for a recent survey. Recently, Pelinovsky--Shimabukuro \cite{PelinovskyShimabukuro2018} proved global well-\-po\-sed\-ness in $H^{1,1}(\mathbb{R}) \cap H^2(\mathbb{R})$ without $L^2$-smallness condition, but assumptions on the Kaup--Newell spectral problem; see also \cite{PelinovskySaalmannShimabukuro2017,Saalmann2017}.

There are major technical difficulties to apply inverse scattering techniques in unweighted Sobolev spaces, e.g., on the line the decay of the data is typically insufficient for classical arguments. For the nonlinear Schr\"odinger equation on the line, Koch--Tataru \cite{Koch2018} could use the transmission coefficient to obtain almost conserved $H^s$-energies for all $s > -\frac{1}{2}$. Killip--Vi\c{s}an--Zhang \cite{KillipVisanZhang2018} pointed out a power series representation for the determinant
\begin{equation*}
\log {\det}
\big(
\begin{bmatrix}
(-\partial +\tilde{\kappa})^{-1} & 0 \\
0 & (- \partial - \tilde{\kappa})^{-1}
\end{bmatrix}
\begin{bmatrix}
-\partial + \tilde{\kappa} & iq \\
\mp i\bar{q} & -\partial - \tilde{\kappa}
\end{bmatrix}
\big),
\end{equation*}
given by
\begin{equation*}
\sum_{l=1}^{\infty} \frac{(\mp 1)^{l-1}}{l} \operatorname{tr}\left\{\left[(\tilde{\kappa}-\partial)^{-1 / 2} q(\tilde{\kappa}+\partial)^{-1} \bar{q}(\tilde{\kappa}-\partial)^{-1 / 2}\right]^{l}\right\},
\end{equation*}
which works in either geometry. Killip \emph{et al.} \cite{KillipVisanZhang2018} showed that it is conserved for NLS and mKdV by term-by-term differentiation. This led to low regularity conservation laws and corresponding a priori estimates in either geometry. Talbut \cite{Talbut2018} used the same approach to show low regularity conservation laws for the Benjamin-Ono equation.

Motivated by these results, we show that the determinant 
\begin{equation*}
\log {\det}
\big(
\begin{bmatrix}
(\partial + i\kappa^2)^{-1} & 0 \\
0 & (\partial - i\kappa^2)^{-1}
\end{bmatrix}
\begin{bmatrix}
\partial + i\kappa^2 & - \kappa q \\
\kappa \bar{q} & \partial - i\kappa^2
\end{bmatrix}
\big),
\end{equation*}
given by
\begin{equation}\label{eq:PerturbationDeterminant}
\sum_{l=1}^{\infty} \frac{(- 1)^{l}i^{l+1}\tilde{\kappa}^{l}}{l} \operatorname{tr}\left\{\left[(\partial - \tilde{\kappa})^{-1} q(\partial+ \tilde{\kappa})^{-1} \bar{q} \right]^{l}\right\},
\end{equation}
where we formally set $\tilde{\kappa} = - i \kappa^2$ (we drop the tilde later on), is conserved for solutions of \eqref{eq:dNLS}. This yields the following theorem on the growth of Besov norms:
\begin{theorem}
\label{thm:APrioriEstimates}
Let $q \in C^\infty(\mathbb{R};\mathcal{S}(\mathbb{K}))$ be a smooth solution to \eqref{eq:dNLS}. For any $0<s<1/2$, $r \in [1,\infty]$, there is $c=c(s,r)<1$ such that
\begin{equation}
\label{eq:APrioriBesovEstimate}
\Vert q(t) \Vert_{B^s_{r,2}} \lesssim \Vert q(0) \Vert_{B^s_{r,2}}
\end{equation}
provided that $\Vert q(0) \Vert_2 \leq c$.
\end{theorem}
\begin{remark}
We focus on regularities for which global results were previously unknown. It appears feasible to cover higher regularities following \cite[Section~3]{KillipVisanZhang2018}.
\end{remark}
In follow-up works to \cite{KillipVisanZhang2018}, Killip--Vi\c{s}an showed sharp global well-\-po\-sed\-ness for the KdV equation \cite{KillipVisan2019} and later on with Bringmann for the fifth order KdV equation \cite{BringmannKillipVisan2019}. Sharp global well-posedness for NLS and mKdV on the real line was shown by Harrop-Griffiths--Killip--Vi\c{s}an \cite{HarropGriffithsKillipVisan2020}. In the first version of the present article (07/2020) we raised the question whether \eqref{eq:dNLS} is within the thrust of these works.

And indeed, in \cite{KillipNtekoumeVisan2021} (01/2021) Killip--Ntekoume--Vi\c{s}an showed global well-\-po\-sed\-ness of \eqref{eq:dNLS} for $q_0 \in H^s(\mathbb{K})$, $\frac{1}{6} \leq s < \frac{1}{2}$ with $\| q_0 \|^2_{L^2} < 4 \pi$. Also, in 12/2020, Tang--Xu \cite{TangXu2020} pointed out an underlying microscopic conservation law on the real line, which paralleled the results in \cite{HarropGriffithsKillipVisan2020} for mKdV and NLS, but did not prove well-posedness. On the real line, Bahouri--Perelman \cite{BahouriPerelman2020} (12/2020) showed global well-posedness in $H^{\frac{1}{2}}(\mathbb{R})$ without smallness assumptions on the $L^2$-norm, relying on profile decomposition, and also crucially on complete integrability. Moreover, Isom--Mantzavinos--Stefanov \cite{IsomMantzavinosStefanov2020} (12/2020) showed that Sobolev norms $H^s(\mathbb{T})$, $s>1$, of solutions to \eqref{eq:dNLS} are growing polynomially by using nonlinear smoothing and not relying on complete integrability.
\vspace{0.5cm}

\textit{Outline of the paper.}
In Section \ref{section:Preliminaries} we introduce notations and recall basic facts about operator traces. In Section \ref{section:PerturbationDeterminant} we show that \eqref{eq:PerturbationDeterminant} is a conserved quantity. In Section \ref{section:ConservationSobolevNorms} we derive low regularity conservation laws, yielding a priori estimates and finishing the proof of the main result.

\section{Preliminaries}
\label{section:Preliminaries}
We start by introducing notations and recalling basic facts about trace class operators and Schatten norms, which will be used in the following.

On the real line, the Fourier transform is defined as
\begin{equation*}
\hat{f}(\xi)= \frac{1}{\sqrt{2 \pi}} \int_{\mathbb{R}} e^{- ix \xi} f(x) dx, \quad f(x) = \frac{1}{\sqrt{2 \pi}} \int_{\mathbb{R}} e^{ix \xi} \hat{f}(\xi) d\xi.
\end{equation*}
The scalar product on $L^2(\mathbb{R})$ is denoted by
\begin{equation*}
    \langle f, g \rangle = \int_{\mathbb{R}} f(x) \overline{g(x)} dx.
\end{equation*}

We will also work on the rescaled torus, for which we use the conventions from \cite{OhWang2020}. Given $\lambda \geq 1$, let $\mathbb{T}_\lambda=\mathbb{R}/(2 \pi \lambda \mathbb{Z})$. The scalar product on $L^2(\mathbb{T}_\lambda)$ is given by
\begin{equation*}
\langle f, g \rangle = \int_{\mathbb{T}_\lambda} f(x) \overline{g(x)} dx = \int_{0}^{2 \pi \lambda} f(x) \overline{g(x)} dx.    
\end{equation*}
We set
\begin{equation*}
\hat{f}(\xi) = \frac{1}{\sqrt{2 \pi}} \int_0^{2 \pi \lambda} f(x) e^{-i x \xi} dx \text{ and } f(x) = \frac{1}{\sqrt{2 \pi}\lambda} \sum_{\xi \in \mathbb{Z}/\lambda} \hat{f}(\xi) e^{i x \xi}
\end{equation*}
for $f \in L^1(\mathbb{T}_\lambda,\mathbb{C})$, where $\xi \in \mathbb{Z}_\lambda = \lambda^{-1} \mathbb{Z}$. The guideline for the conventions is that Plancherel's theorem remains true:
\begin{equation*}
\Vert f \Vert_{L^2(\mathbb{T}_\lambda)} = \Vert \hat{f} \Vert_{L^2(\mathbb{Z}_\lambda,(d\xi)_\lambda)},
\end{equation*}
where $(d\xi)_\lambda$ denotes the normalized counting measure on $\mathbb{Z}_\lambda$:
\begin{equation*}
\int_{\mathbb{Z}_\lambda} f(\xi) (d\xi)_\lambda = \frac{1}{\lambda} \sum_{\xi \in \mathbb{Z}_\lambda} f(\xi).
\end{equation*}
For further basic Fourier analysis identities on $\mathbb{T}_\lambda$, we refer to \cite[Section~2]{CollianderKeelStaffilaniTakaokaTao2002}. We turn to the definition of $L^2$-based Sobolev norms: For $s \in \mathbb{R}$, $u \in \mathcal{S}(\mathbb{R})$ we consider
\begin{equation*}
\| f \|_{H^s(\mathbb{R})} = \left( \int_{\mathbb{R}} (1+|\xi|^2)^{s} |\hat{f}(\xi)|^2 d\xi \right)^{\frac{1}{2}},
\end{equation*}
or for $f \in \mathcal{S}(\mathbb{T}_\lambda) = C^\infty(\mathbb{T}_\lambda)$
\begin{equation*}
\| f \|_{H^s(\mathbb{T}_\lambda)} = \left( \int_{\mathbb{Z}_\lambda} (1+|\xi|^2)^s |\hat{f}(\xi)|^2 (d\xi)_\lambda \right)^{\frac{1}{2}}.
\end{equation*}
For the definition of Besov norms, we consider a smooth partition of unity of the real line: Let $\beta_1: \mathbb{R} \to [0,1]$ denote a radially decreasing function $\beta_1(\xi) = 1$ for $\xi \in [-1,1]$ and $\text{supp } \beta_1 \subseteq [-2,2]$. For $N \in 2^{\mathbb{N}}$ let $\beta_N(\xi) = \beta_1(\xi/N) - \beta_1(\xi/(N/2))$, and let $P_N$ denote the Fourier multiplier on $\mathbb{R}$ or $\mathbb{T}_\lambda$:
\begin{equation*}
(P_N f) \widehat (\xi) = \beta_N(\xi) \hat{f}(\xi).
\end{equation*}
We define the Besov norm of $f \in \mathcal{S}(\mathbb{R})$ or $f \in \mathcal{S}(\mathbb{T}_\lambda)$ for $1 \leq r < \infty$, $s \geq 0$  by
\begin{equation*}
\| f \|_{B^s_{r,2}} = \big( \sum_{N \in 2^{\mathbb{N}_0}} N^{rs} \| P_N f \|^r_{L^2} \big)^{\frac{1}{r}}
\end{equation*}
and with the usual modification for $r = \infty$.\\
For $\lambda \in 2^{\mathbb{N}_0}$, let $f_\lambda(x) = \lambda^{-\frac{1}{2}} f(\lambda^{-1} x)$. We record the following scaling of the Besov norms:
\begin{equation}
\label{eq:BesovScaling}
\lambda^{-s} \| f \|_{B^s_{r,2}} \lesssim \| f_\lambda \|_{B^s_{r,2}} \lesssim \| f \|_{L^2} + \lambda^{-s} \| f \|_{B^s_{r,2}},
\end{equation}
which follows from 
\begin{align*}
\| P_1 f_\lambda \|_{L^2} \leq \| f \|_{L^2}, \quad \| P_N f_\lambda \|_{L^2} = \| P_{\lambda N} f \|_{L^2} \qquad (N \in 2^{\mathbb{N}}).
\end{align*}

For $\kappa > 0$, $k \in \mathbb{Z}$, the mappings $( \partial \pm \kappa)^{k}: \mathcal{S}(\mathbb{K}) \to \mathcal{S}(\mathbb{K})$, $\mathbb{K} \in \{ \mathbb{R}, \mathbb{T}_\lambda \}$, are defined as Fourier multipliers:
\begin{equation*}
(( \partial \pm \kappa)^{k} f ) \widehat (\xi) = (i \xi \pm \kappa)^k \hat{f}(\xi).
\end{equation*}
Bounds in $L^2$-based Sobolev spaces are immediate from Plancherel's theorem.
We denote $R_{\pm} = (\partial \pm \kappa)^{-1}$, which have the following kernels on the real line:
\begin{equation}
\label{eq:KernelsR}
k_+(\kappa,x,y) = 
\begin{cases}
e^{-\kappa(x-y)} & \text{ if } x > y, \\
0 & \text{ else},
\end{cases}
\quad k_-(\kappa,x,y) = 
\begin{cases}
 -e^{\kappa(x-y)} & \text{ if } x <y, \\
0 & \text{ else}.
\end{cases}
\end{equation}
On the circle, by the Poisson summation formula (cf. \cite[Lemma~3.3]{OhWang2020}) we find the kernels of $R_{\pm}$ to be
\begin{equation}
\label{eq:KernelsT}
k^\lambda_{-}(\kappa,x,y) = - \frac{e^{\kappa((x-y) - 2\pi \lambda \lceil \frac{x-y}{2 \pi \lambda} \rceil)}}{1 - e^{-2 \pi \lambda \kappa}}, \quad k^\lambda_+(\kappa,x,y) = \frac{e^{\kappa((y-x) - 2\pi \lambda \lceil \frac{y-x}{2 \pi \lambda} \rceil)}}{1 - e^{-2 \pi \lambda \kappa}},
\end{equation}
where $\lceil \cdot \rceil: \mathbb{R} \to \mathbb{Z}$ denotes the ceiling function given by $\lceil x \rceil = \min \{ k \in \mathbb{Z} : k \geq x \}$. We note the following identity:
\begin{equation}
\label{eq:K-2Relation}
    k_-^\lambda(\kappa,x,y)^2 = \frac{1+e^{-2 \pi \lambda \kappa}}{1- e^{-2 \pi \lambda \kappa}} k_-^\lambda(2\kappa,x,y).
\end{equation}

Let $H$ be a Hilbert space. For an introduction to the following concepts regarding compact operators on $H$, we refer to \cite[Chapter~3]{Simon2015}. By $\mathfrak{I}_p(H)$ we denote the Schatten class of compact operators with $\ell^p$-summable singular values. From the embedding properties of the sequence space $\ell^p$ we see the embeddings $\mathfrak{I}_p \subset \mathfrak{I}_q$ when $p < q$. The space $\mathfrak{I}_\infty$ is the space of compact operators. By Hilbert Schmidt (HS) operators we refer to elements of $\mathfrak{I}_2(H)$; elements of $\mathfrak{I}_1(H)$ are referred to as of trace-class. These spaces are *-ideals in the space of bounded linear operators on $H$, and the following estimate holds:
\begin{equation*}
    \|AB\|_{\mathfrak{I}_p} + \|BA\|_{\mathfrak{I}_p} \leq \|A\|_{\mathfrak{I}_p} \|B\|_{H\to H}.
\end{equation*}

Next, we explain how to compute the operator trace of products of two Hilbert Schmidt operators, which suffices for the present context. Firstly, recall that an operator $A:L^2(\mathbb{K}) \to L^2(\mathbb{K})$ is HS if and only if it admits a kernel $K(x,y) \in L^2(\mathbb{K} \times \mathbb{K})$ and
\begin{equation*}
\| A \|_{L^2(\mathbb{K}) \to L^2(\mathbb{K})} \leq \| A \|_{\mathfrak{I}_2(\mathbb{K})} = \iint_{\mathbb{K} \times \mathbb{K}} |K(x,y)|^2 dx dy.
\end{equation*}
Furthermore, products of two or more HS operators $A_i:L^2(\mathbb{K}) \to L^2(\mathbb{K})$ with kernels $K_i \in L^2(\mathbb{K} \times \mathbb{K})$ are of trace class, and the operator trace is computed as
\begin{equation*}
    \text{tr}(A_1 \ldots A_n) = \int_{\mathbb{K}} \ldots \int_{\mathbb{K}} K_1(x_1,x_2) K_2(x_2,x_3) \ldots K_n(x_n,x_1) dx_n \ldots dx_1.
\end{equation*}
By Fubini's theorem, this allows to cycle the trace:
\begin{equation*}
    \text{tr}(A_1 \ldots A_n) = \text{tr}(A_2 \ldots A_n A_1).
\end{equation*}
Moreover, for a trace class operator $A$ and a bounded operator $B$ on $L^2(\mathbb{K})$, the operators $AB$ and $BA$ are of trace class, and we can cycle the trace:
\begin{equation*}
    \text{tr}(AB) = \text{tr}(BA).
\end{equation*}
We will also need the H\"{o}lder-like estimate for Schatten norms
\begin{equation*}
    \|AB\|_{\mathfrak{I}_p} \leq \|A\|_{\mathfrak{I}_r}\|B\|_{\mathfrak{I}_s}
\end{equation*}
provided that $1/p = 1/r + 1/s$, $1 \leq p,r,s \leq \infty$.

\section{The perturbation determinant}
\label{section:PerturbationDeterminant}
In this section we show conservation of the perturbation determinant 
\begin{equation}
\label{eq:RealPerturbationDeterminant}
\alpha(\kappa;q) = \mathrm{Re} \sum_{l \geq 1} \frac{ (-i)^{l+1} \kappa^l}{l} \text{tr} ((\partial - \kappa)^{-1} q (\partial + \kappa)^{-1} \bar{q})^l) = \sum_{l \geq 1} \alpha_l
\end{equation}
through term-by-term differentiation.
For the first term we note the following:
\begin{lemma}
\label{lem:R-R+ScalarProduct}
    The following identities hold for $f,g \in \mathcal{S}$:
    \begin{equation}
        \operatorname{tr}((\kappa - \partial)^{-1}f(\kappa + \partial)^{-1}g) = \begin{cases}
            \langle (2\kappa - \partial)^{-1}f,\bar{g}\rangle,\quad \text{if } \mathbb{K} = \mathbb{R},\\
            \frac{1+e^{-2 \pi \lambda \kappa}}{1-e^{-2 \pi \lambda \kappa}} \langle (2\kappa - \partial)^{-1}f,\bar{g}\rangle, \quad \text{if } \mathbb{K} = \mathbb{T}_\lambda.
        \end{cases}
    \end{equation}
\end{lemma}
\begin{proof}
    We begin with the line case. Using the explicit kernels \eqref{eq:KernelsR}, we find
    \begin{align*}
        \operatorname{tr}((\kappa - \partial)^{-1}f(\kappa + \partial)^{-1}g) &= \iint_{\mathbb{R}^2 \cap \{x<y\}} e^{2\kappa(x-y)}f(y)g(x) \, dx dy \\
        &= \langle (2\kappa - \partial)^{-1}f, \bar{g}\rangle,
    \end{align*}
    the last line using the $L^2(\mathbb{R})$ scalar product. In the circle case, \eqref{eq:K-2Relation} yields
    \begin{align*}
        \operatorname{tr}((\kappa - \partial)^{-1}f(\kappa + \partial)^{-1}g) &= \iint_{\mathbb{T}^2_\lambda} k_-^\lambda(\kappa,x,y)f(y)k_+^\lambda(\kappa,y,x)g(x) \, dx dy\\
        &= \frac{1+e^{-2 \pi \lambda \kappa}}{1-e^{-2 \pi \lambda \kappa}}\iint_{\mathbb{T}^2_\lambda} k_-^\lambda(2\kappa,x,y)f(y)g(x) \, dx dy\\
        &= \frac{1+e^{-2 \pi \lambda \kappa}}{1-e^{-2 \pi \lambda \kappa}} \langle (2\kappa - \partial)^{-1}f,\overline{g}\rangle.
    \end{align*}
\end{proof}

 To ensure that we can differentiate $\alpha$ term by term, we show the following trace estimates leading to geometric convergence:
\begin{lemma}
\label{lem:TraceEstimates}
Let $q \in \mathcal{S}(\mathbb{R})$ or $q \in \mathcal{S}(\mathbb{T}_\lambda)$, $\lambda \geq 1$, $l \geq 2$, and $\kappa > 0$. Then, we find the following estimates to hold:
\begin{equation}
\label{eq:TraceEstimates}
\begin{split}
    \left| \operatorname{tr} \{ [ \kappa (\kappa - \partial)^{-1} q (\kappa + \partial)^{-1} \bar{q}]^l \} \right| \lesssim
 \begin{cases}
(\kappa^{-s} \| q \|_{H^s})^{2l} \quad (0 \leq s < 1/4), \\
(\kappa^{-1/4} \| q \|_{H^s} )^{2l} \quad (s \geq 1/4)
\end{cases}
\end{split}
\end{equation}

Hence, for $\Vert q \Vert_{L^2} \leq c \ll 1$ small enough, or $\kappa \gg \Vert q \Vert_{H^{s}}^{1/a}$ for $s>0$ and $a=\min(1/4,s)$, $\alpha$ defined in \eqref{eq:RealPerturbationDeterminant} converges geometrically.
\end{lemma}
\begin{proof} 
The $L^2$-estimate in \eqref{eq:TraceEstimates} follows by
 \begin{equation}
 \label{eq:HSEstimateII}
\| R_{\pm} q \|^2_{\mathfrak{I}_2} \sim \frac{1}{\kappa} \| q \|^2_{L^2}, 
 \end{equation}
 To show the above display, we note that the kernel is given by $K(x,y) = k_{\pm}(x,y) q(y)$, and we compute by Fubini and Plancherel's theorem
\begin{equation*}
\begin{split}
\iint |K(x,y)|^2 dx dy &= \int dy |q(y)|^2 \int dx |k_{\pm}(x,y)|^2 \\
 &= \int dy |q(y)|^2 \int \frac{d\xi}{\kappa^2 + \xi^2} \sim \frac{\| q \|^2_{L^2}}{\kappa}.
\end{split}
\end{equation*}
Observe that the argument works in either geometry and in the periodic case gives a bound independent of the period length. 

For the $H^s$-part, set $A = \kappa^{1/2} R_{\pm} q$. Firstly, we argue that we can gain powers of $\kappa$ by estimating $q$ in $H^s$-norms. Note that 
\begin{align}
\label{eq:EstimatesA}
\Vert A \Vert_{\mathfrak{I}_2} \lesssim \Vert q \Vert_2, \qquad \Vert A \Vert_{\mathfrak{I}_\infty} \lesssim \kappa^{-1/2} \Vert q \Vert_{\infty}.
\end{align}
The first estimate follows from \eqref{eq:HSEstimateII} and the second follows from viewing $q$ as a multiplication operator in $L^2$ and the bound $\| R_{\pm} \|_{L^2 \to L^2} \lesssim \kappa^{-1}$, which is immediate from Plancherel's theorem. Interpolating the estimates in \eqref{eq:EstimatesA} by viewing $A$ as a bounded operator from $L^p \to \mathfrak{I}_p$ (cf. \cite[Proposition~I.1]{BenyaminiLindenstrauss2000}) and using Sobolev embedding, we find
\begin{equation*}
\Vert A \Vert_{\mathfrak{I}_p} \lesssim \kappa^{-1/2+1/p} \Vert q \Vert_p \lesssim \kappa^{-s} \Vert q \Vert_{H^s} \text{ for } s = 1/2-1/p, \quad 2 \leq p < \infty.
\end{equation*}

Let $0<s'<1/4$ in the following and set $s'=\frac{1}{2} - \frac{1}{p_{s'}}$. By H\"older's inequality and embeddings for Schatten spaces, we find
\begin{equation*}
    |\alpha_2| \leq |\text{tr} ((A \bar{A})^2 ) | \leq \Vert A \Vert^4_{\mathfrak{I}_4} \leq \Vert A \Vert^4_{\mathfrak{I}_{p_{s^\prime}}} \lesssim \kappa^{-4s^\prime} \Vert q \Vert_{H^{s^\prime}}^4.
\end{equation*}
Similarly for the higher order terms $l \geq 3$, we find
\begin{equation*}
|\text{tr} ((A \bar{A})^{l}) | \leq \Vert A \Vert^{2l}_{\mathfrak{I}_{2l}} \leq \Vert A \Vert_{\mathfrak{I}_{p_{s'}}}^{2l} 
\lesssim \kappa^{-2ls'} \Vert q \Vert_{H^{s'}}^{2l}.
\end{equation*}
This is the $H^s$-estimate in \eqref{eq:TraceEstimates} for $0<s<1/2$. Lastly, note that this implies that the series \eqref{eq:RealPerturbationDeterminant} converges for $q \in H^s$, $s>0$ by choosing $\kappa \gg \Vert q \Vert_{H^{s}}^{1/a}$ for $a=\min(1/4,s)$ as claimed.
 \end{proof}

Next, we show that $\alpha$ is conserved by term-by-term differentiation.
\begin{proposition}
\label{prop:ConservationPerturbationDeterminant}
Let $q \in C^\infty(\mathbb{R};\mathcal{S})$ be a smooth global solution to \eqref{eq:dNLS} with \\
 $\Vert q(0) \Vert_2 \leq c \ll 1$. Then,
\begin{equation*}
\frac{d}{dt} \alpha(\kappa;q) = 0.
\end{equation*}
\end{proposition}
\begin{remark}
\label{rem:HsConvergence}
By Lemma \ref{lem:TraceEstimates} $\alpha(\kappa;q)$ converges \emph{without smallness assumption} on the $L^2$-norm, but provided that $\kappa$ is sufficiently large. However, we are not able to show bounds for the $B_{r,2}^s$-norm without $L^2$-smallness assumption. 
\end{remark}
\begin{proof}
In the following we omit taking the real part in \eqref{eq:PerturbationDeterminant} and will thus show that both real and imaginary part are conserved. Consider
\begin{equation*}
\sum_{l=1}^\infty \frac{(-i)^{l+1}  \kappa^l }{l} \text{tr} (\underbrace{( \partial - \kappa)^{-1}}_{R_-} q \underbrace{(\partial + \kappa)^{-1}}_{R_+} \bar{q})^l = \sum_{l \geq 1} \alpha_l.
\end{equation*}
We note similar to the considerations from \cite[Section~4]{KillipVisanZhang2018}:
\begin{equation}
\label{eq:q^2qOperatorIdentity}
\begin{split}
(|q|^2 q)_x &= (\partial- \kappa)(|q|^2 q) -(|q|^2 q) (\partial + \kappa) + 2 \kappa (|q|^2 q) \\
(|q|^2 \bar{q})_x &= (\partial + \kappa)(|q|^2 \bar{q}) - (|q|^2 \bar{q})(\partial - \kappa) - 2 \kappa (|q|^2 \bar{q}),
\end{split}
\end{equation}
and furthermore,
\begin{equation}
\label{eq:qxxOperatorIdentity}
\begin{split}
q_{xx} &= q (\partial^2 - 2 \kappa \partial - \kappa^2) + (\partial^2 + 2\kappa \partial - \kappa^2) q + 2(\kappa - \partial) q (\kappa+\partial), \\
\bar{q}_{xx} &= (\partial^2 - 2 \kappa \partial - \kappa^2) \bar{q} + \bar{q} (\partial^2 + 2 \kappa \partial - \kappa^2) + 2(\kappa + \partial)\bar{q} (\kappa - \partial).
\end{split}
\end{equation}
Differentiating term-by-term, we find two terms $\frac{d}{dt} \alpha_l = A_l + B_l$, which are given by
\begin{align*}
A_l &= -(-i)^{l+1} \kappa^l \text{tr} ((R_- q R_+ \bar{q})^{l-1} [R_- (|q|^2 q)_x R_+ \bar{q} + R_- q R_+ (|q|^2 \bar{q})_x] \\
B_l &= (-i)^{l+1} \kappa^l \text{tr} ((R_- q R_+ \bar{q})^{l-1} [R_- i q_{xx} R_+ \bar{q} - i R_- q R_+ \bar{q}_{xx} ]).
\end{align*}
We show that $A_l + B_{l+1} = 0$ by substituting \eqref{eq:q^2qOperatorIdentity} and \eqref{eq:qxxOperatorIdentity}. However, with the substitutions introducing differential operators, we have to check that the single terms are well-defined and cycling the trace is admissible. Strictly speaking, already writing $A_l$ and $B_l$ as in the above display requires cycling the trace. Since $R_- q R_+ \bar{q}$, $R_-(|q|^2 q)_x R_+ \bar{q}$, and $R_- q_{xx} R_+ \bar{q}$ are of trace class, this is not an issue. We shall prove $A_l + B_{l+1} = 0$ for $l \geq 2$ by substitution and handle the terms $A_1$, $B_1$, and $B_2$ directly.

For $A_l$ we find after substitution of \eqref{eq:q^2qOperatorIdentity}:
\begin{equation}
\label{eq:Al}
\begin{split}
A_l = - (-i)^{l+1} \kappa^l \text{tr} (&(R_- q R_+ \bar{q})^{l-1} [|q|^2 q R_+ \bar{q} - R_- |q|^4 + 2 \kappa R_- |q|^2 q R_+ \bar{q}  \\
+ &\; R_- |q|^4 - R_- q R_+ |q|^2 \bar{q} R_-^{-1} - 2 \kappa R_- q R_+ |q|^2 \bar{q} ]).
\end{split}
\end{equation}
With $R_- q R_+ \bar{q}$ being trace class, for $l \geq 2$ it is enough to check boundedness of the remaining six factors. With $R_-$ and multiplication by $q$ or $\bar{q}$ a bounded operator as $q \in \mathcal{S}$, it only remains to check boundedness of the fifth factor: $R_- q R_+ |q|^2 \bar{q} R_-^{-1}$. Now $R_-^{-1}: L^2 \to H^{-1}$ is bounded and so it is enough to see that multiplication with $|q|^2 \bar{q}$ is bounded in $H^{-1}$ because $R_+$ is a bounded operator $H^{-1} \to L^2$. This follows from the estimate $\|fg\|_{H^{-1}} \lesssim \|f\|_{H^1}\|g\|_{H^{-1}}$, which is immediate by duality and the algebra property of $H^1(\mathbb{K})$. Hence, we can consider the traces of the single terms, and the second cancels the fourth term. For the fifth term, we compute by cycling the trace
\begin{equation*}
\begin{split}
&\; \text{tr} (R_- q R_+ \bar{q} (R_- q R_+ \bar{q})^{l-2} R_- q R_+ |q|^2 \bar{q} R_{-}^{-1} ) \\
&= \text{tr} ((R_- q R_+ \bar{q})^{l-2} R_- q R_+ |q|^2 \bar{q} R_{-}^{-1} R_- q R_+ \bar{q}) \\
&= \text{tr} ((R_- q R_+ \bar{q})^{l-2} R_- q R_+ |q|^4 R_+ \bar{q} ),
\end{split}
\end{equation*}
which thus cancels the first term.
Hence, for $l \geq 2$, we have proved
\begin{equation}
\label{eq:Al2}
A_l = - (-i)^{l+1} \kappa^l \text{tr} ((R_- q R_+ \bar{q})^{l-1} [ 2 \kappa R_- |q|^2 q R_+ \bar{q} - 2 \kappa R_- q R_+ |q|^2 \bar{q} ]).
\end{equation}
For $B_{l+1}$, $l \geq 1$, we find after substitution of \eqref{eq:qxxOperatorIdentity}:
\begin{equation}
\begin{split}
\label{eq:Bl+1}
B_{l+1} &= (-i)^{l+1} \, \kappa^{l+1} \, \text{tr} ((R_- q R_+ \bar{q})^l [ R_- q(\partial^2 - 2\kappa \partial - \kappa^2) R_+ \bar{q} \\
&\quad + R_- (\partial^2 + 2 \kappa \partial - \kappa^2) q R_+ \bar{q} - 2 |q|^2 -  R_- q R_+ (\partial^2 - 2\kappa \partial -\kappa^2) \bar{q} \\
 &\quad - R_- q R_+ \bar{q} (\partial^2 + 2 \kappa \partial - \kappa^2) - 2R_- |q|^2 (\kappa - \partial)]).
\end{split}
\end{equation}
We have to verify that the traces of the single terms are well-defined, for which it is again enough to see the boundedness of the six factors with $R_- q R_+ \bar{q}$ being trace class. This follows similarly to the above. Consider e.g. the first term $R_- q (\partial^2 - 2 \kappa \partial + \kappa^2) R_+ \bar{q}$. With $(\partial^2 - 2 \kappa \partial + \kappa^2) R_+: L^2 \to H^{-1}$ bounded and multiplication with $q$ or $\bar{q}$ bounded in $H^s$, $s \in \mathbb{R}$, we find that $q (\partial^2 - 2 \kappa \partial + \kappa^2) R_+ \bar{q}: L^2 \to H^{-1}$ is bounded. Composition with $R_-$ yields $L^2$-boundedness.

Next, observe that the first and fourth term cancel because constant coefficient differential operators are commuting. This also implies cancelling of the second and fifth term, after additionally cycling the trace:
\begin{equation*}
\begin{split}
&\quad \text{tr} ((R_- q R_+ \bar{q})(R_- q R_+ \bar{q})^{l-1} R_- q R_+ \bar{q} (\partial^2 + 2 \kappa \partial - \kappa^2)] \\
&= \text{tr} ((R_- q R_+ \bar{q})^{l} (\partial^2 + 2 \kappa \partial - \kappa^2) R_- q R_+ \bar{q}).
\end{split}
\end{equation*}
To summarize, we have found
\begin{equation}
\label{eq:Bl+12}
B_{l+1} = (-i)^{l+1} \, 2 \kappa^{l+1} \, \text{tr} ((R_- q R_+ \bar{q})^l [- |q|^2 + R_- |q|^2 R_-^{-1}]).
\end{equation}
The first term in \eqref{eq:Bl+12} is cancelled by the second term of \eqref{eq:Al2}, and the second term in \eqref{eq:Bl+12} is cancelled by the first term of \eqref{eq:Al2} after additionally cycling the trace.

It remains to prove $B_1 = 0$ and $A_1 + B_2 = 0$. The first claim was already shown in \cite[Eq.~(51)]{KillipVisanZhang2018} and follows from similar considerations as below.

We turn to the second claim: The substitution \eqref{eq:q^2qOperatorIdentity} in $A_1$ cannot easily be justified\footnote{We thank the referee for pointing this out.}; for $B_2$ \eqref{eq:Bl+12} remains correct. Hence, we resort to doing the integration by parts in $A_1$ directly. 
Recall
\begin{equation*}
A_1 = \kappa \text{tr} [ R_- (|q|^2 q)_x R_+ \bar{q} + R_- q R_+ (|q|^2 \bar{q})_x ].
\end{equation*}
By Lemma \ref{lem:R-R+ScalarProduct}, we find
\begin{equation*}
\begin{split}
    \text{tr} (R_- (|q|^2 q)_x R_+ \overline{q}) &= - \langle (2\kappa - \partial)^{-1} (|q|^2 q)_x, q \rangle \\
    &= - \langle (2\kappa - \partial)^{-1} [(\partial - 2 \kappa) + 2\kappa] (|q|^2 q), q \rangle \\
    &=  \langle |q|^2 q, q \rangle - 2 \kappa \langle (2\kappa - \partial)^{-1} (|q|^2 q), q \rangle \\
    &= \int |q|^4 dx + 2 \kappa \text{tr} (R_- |q|^2 q R_+ \bar{q}).
    \end{split}
\end{equation*}
In a similar spirit, we compute
\begin{equation*}
\begin{split}
    \text{tr} (R_- q R_+ (|q|^2 \bar{q})_x) &= - \langle (2\kappa - \partial)^{-1} q, (|q|^2 q)_x \rangle \\
    &=  \langle (2\kappa - \partial)^{-1} q_x, |q|^2 q \rangle \\
    &=  \langle (2\kappa - \partial)^{-1} (\partial - 2 \kappa + 2\kappa)q, |q|^2 q \rangle \\
    &= - \langle q, |q|^2 q \rangle + 2 \kappa \langle (2\kappa - \partial)^{-1} q, |q|^2 q \rangle \\
    &= - \int |q|^4 dx - 2\kappa \text{tr} (R_- q R_+ |q|^2 \bar{q}).
    \end{split}
\end{equation*}
With the $L^4$-norms cancelling, we conclude
\begin{equation*}
A_1 = 2 \kappa^2 ( \text{tr} (R_- |q|^2 q R_+ \bar{q} - R_- q R_+ |q|^2 \bar{q})) = -B_2.    
\end{equation*}
The proof is complete.
\end{proof}

\section{Conservation of Besov norms with positive regularity index}
\label{section:ConservationSobolevNorms}
In the following, we want to construct Besov norms from the leading term of $\alpha(q;\kappa)$. Set
\begin{equation*}
w(\xi,\kappa) = \frac{\kappa^2}{\xi^2 + 4 \kappa^2} - \frac{(\kappa/2)^2}{\xi^2 + \kappa^2} = \frac{3 \kappa^2 \xi^2}{4(\xi^2+ \kappa^2)(\xi^2+ 4 \kappa^2)}
\end{equation*}
and
\begin{equation}
\label{eq:BesovSubstitute}
\Vert f \Vert_{Z_r^s} = \big( \sum_{N \in 2^\mathbb{N}} N^{rs} \langle f, w(-i \partial_x, N) f \rangle^{r/2} \big)^{1/r}.
\end{equation}
The $Z_r^s$--norm consists of homogeneous components, which can be linked to the perturbation determinant. We will use the identities \cite[Eq.~(40),~(55)]{KillipVisanZhang2018}:
\begin{align}
\Vert f \Vert_{B^s_{r,2}} &\lesssim_s \Vert f \Vert_{H^{-1}} + \Vert f \Vert_{Z^s_{r}}, \label{eq:BesovZsEmbedding}\\
\Vert f \Vert_{Z^s_{r}} &\lesssim \Vert f \Vert_{B_{r,2}^s}.
\label{eq:ZsBesovEmbedding}
\end{align}
Consequently, it suffices to control the $Z_r^s$-norm to infer about the Besov norms.

\begin{remark}
	In the $Z^s_r$-quantities introduced in \cite{KillipVisanZhang2018}, there is an additional parameter $\kappa_0$. One might hope that this flexibility helps to obtain a result for arbitrary initial data. However, $\kappa_0$ enters with a positive exponent into the estimates. This reflects indeed the relation of $\kappa_0$ with rescaling and the $L^2$-criticality of \eqref{eq:dNLS}. To keep things simple, we choose $\kappa_0 = 1$.
\end{remark}

\subsection{The line case}
To analyze the growth of the $Z_r^s$-norm, we link the multiplier from above with the first term of $\alpha$. We recall the following identity on the real line, which is immediate from Lemma \ref{lem:R-R+ScalarProduct}:
\begin{corollary}
For $\kappa > 0$ and $q \in \mathcal{S}$, we find
\begin{equation*}
\mathrm{Re} \big( \kappa \, \mathrm{tr} \, ( (\kappa - \partial)^{-1} q (\kappa + \partial)^{-1} \bar{q} ) \big) = \int \frac{2 \kappa^2 |\hat{q}(\xi)|^2}{\xi^2+ 4 \kappa^2} d\xi.
\end{equation*}
\end{corollary}
This yields
\begin{equation*}
\begin{split}
\langle f, w(-i \partial_x,N) f \rangle &= \int \frac{N^2}{\xi^2 + 4N^2} |\hat{f}(\xi)|^2 d\xi - \int \frac{(N/2)^2}{\xi^2 + N^2} |\hat{f}(\xi)|^2 d\xi \\
&= \frac{1}{2} [ \alpha_1(N, f) - \alpha_1(N/2,f)].
\end{split}
\end{equation*}

We can estimate $|\alpha - \alpha_1|$ favorably by \eqref{eq:TraceEstimates}
\begin{equation}
\label{eq:RemainderEstimate}
\big| \sum_{l \geq 2} \alpha_l(\kappa,q(t)) \big| \lesssim \kappa^{-4s^\prime} \Vert q(t) \Vert_{H^{s^\prime}}^4
\end{equation}
provided that $0<s'<1/4$ and $\| q(t)\|_{H^{s'}} \leq d_{s'} \ll 1$. Let $D_{s,r}$ denote the constant such that
\begin{equation}
\label{eq:QuantificationBesovEmbedding}
\| q(t) \|_{B^s_{r,2}} \leq D_{s,r} ( \| q(0) \|_{L^2} + \| q(t) \|_{Z^s_r})
\end{equation}
by \eqref{eq:BesovZsEmbedding} and $L^2$-conservation.

\eqref{eq:RemainderEstimate} gives by the embedding $B_{r,2}^{s} \hookrightarrow H^{s^\prime}$ for $s > s^\prime$ and $r \in [1,\infty]$
\begin{align*}
&\qquad \quad \langle q(t), w(-i \partial_x, N) q(t) \rangle \\
 &\lesssim \langle q(0), w(-i\partial_x, N) q(0) \rangle + N^{-4s^\prime} [ \Vert q(t) \Vert_{H^{s^\prime}}^4 + \Vert q(0) \Vert_{H^{s^\prime}}^4 ] \\ 
&\lesssim \langle q(0), w(-i\partial_x, N) q(0) \rangle + N^{-4s^\prime} [ \Vert q(t) \Vert_{B_r^{s,2}}^4 + \Vert q(0) \Vert_{B_r^{s,2}}^4 ].
\end{align*}
Raising the estimate to the power $r/2$, multiplying with $N^{rs}$, and carrying out the dyadic sums over $N \in 2^{\mathbb{N}_0}$, we find
\begin{align*}
\Vert q(t) \Vert^r_{Z^s_r} &\lesssim \Vert q(0) \Vert^r_{Z^s_r} + [ \Vert q(t) \Vert^{2r}_{B_{r,2}^{s}} + \Vert q(0) \Vert^{2r}_{B_{r,2}^{s}} ] 
\end{align*}
provided that we choose $s'<s<2s'$. This can be satisfied for $0<s<1/2$.\\
By \eqref{eq:ZsBesovEmbedding} and $L^2$-conservation, we arrive at
\begin{equation}
\label{eq:PropagationZsNorm}
\Vert q(t) \Vert_{Z^s_r} \leq C_{r,s} ( \Vert q(0) \Vert_{Z^s_r} + \Vert q(0) \Vert^2_{L^2} + [\Vert q(t) \Vert^2_{Z_r^{s}} + \Vert q(0) \Vert^2_{Z_r^{s}} ])
\end{equation}
with $C_{r,s} \geq 1$ provided that $\| q(t') \|_{B^s_{r,2}}$ is small enough for $t' \in [-t,t]$ such that \eqref{eq:RemainderEstimate} holds by $\| q(t') \|_{H^{s'}} \leq \| q(t') \|_{B^s_{r,2}}$.

\eqref{eq:PropagationZsNorm} can be bootstrapped. Suppose that 
\begin{equation*}
\max( \Vert q(0) \Vert_{Z^s_r}, \Vert q(0) \Vert_{L^2}) \leq \varepsilon \ll 1,
\end{equation*}
where $\varepsilon$ is chosen below as $\varepsilon = \varepsilon(C_{r,s}, D_{r,s})$. We prove that for any $t \in \mathbb{R}$
\begin{equation}
\label{eq:BootstrapAssumption}
\Vert q(t) \Vert_{Z^s_r} \leq 2 C_{r,s} \varepsilon.
\end{equation}
For this purpose, let $I$ denote the maximal interval containing the origin such that \eqref{eq:BootstrapAssumption} holds for any $t \in I$. $I$ is non-empty and closed due to continuity of $\Vert q(t) \Vert_{Z^s_r}$. Furthermore, $I$ is open:
For $t \in I$ \eqref{eq:PropagationZsNorm} yields
\begin{equation*}
\Vert q(t) \Vert_{Z^s_r} \leq C_{r,s} ( \varepsilon + 2\varepsilon^2 + (2 C_{r,s} \varepsilon^2)) \leq (3/2) C_{r,s} \varepsilon
\end{equation*}
by choosing $\varepsilon \leq (8 C_{r,s})^{-1}$, but \eqref{eq:PropagationZsNorm} hinges on veracity of \eqref{eq:RemainderEstimate}. To guarantee this, we choose $\varepsilon$ possibly smaller such that 
\begin{equation*}
 \varepsilon D_{r,s} (1+ 2 C_{r,s}) \leq c_{s'}.
\end{equation*}
By $\| q(t) \|_{H^{s'}} \leq \|q(t) \|_{B^s_{r,2}}$ and \eqref{eq:QuantificationBesovEmbedding}, this additionally shows that \eqref{eq:PropagationZsNorm} is true. We conclude that $I = \mathbb{R}$.

This finishes the proof for initial data with small $L^2$- and $Z^s_r$-norm. The assumption $\Vert q(0) \Vert_{Z^s_r} \leq \varepsilon$ follows for initial data with smaller $L^2$-norm through rescaling $q_0(x) = \lambda^{-\frac{1}{2}} q_0(\lambda^{-1} x)$, $\lambda \in 2^{\mathbb{N}}$ by \eqref{eq:ZsBesovEmbedding} and \eqref{eq:BesovScaling}. The proof of Theorem \ref{thm:APrioriEstimates} is complete in the line case.
\subsection{The circle case}
In this section we discuss the case of periodic boundary conditions.
We shall rescale the circle, too, to accomplish smallness of the homogeneous norms. In \cite{KillipVisanZhang2018}, this was not necessary due to more freedom in the parameter $\kappa$. 
For the leading term in \eqref{eq:PerturbationDeterminant} we find from Lemma \ref{lem:R-R+ScalarProduct} (see also \cite[Lemma~3.3]{OhWang2020}):
\begin{corollary}
\label{lem:LeadingTermCircle}
Let $\kappa \geq 1$ and $\lambda \geq 1$. Then, we have
\begin{equation*}
    \mathrm{Re} \operatorname{tr} \big( \kappa (\kappa - \partial)^{-1} q (\kappa + \partial)^{-1} \bar{q} \big) = \frac{1 + e^{-2 \pi \lambda \kappa}}{1-e^{-2 \pi \lambda \kappa}} \int_{\mathbb{Z}_\lambda} \frac{2 \kappa^2 |\hat{q}(\xi)|^2}{4 \kappa^2 + \xi^2} (d\xi)_\lambda
\end{equation*}
for any smooth function $q$ on $\mathbb{T}_\lambda$.
\end{corollary}

Set $C(\lambda,N) = (1+e^{-2 \pi \lambda N})/(1-e^{-2 \pi \lambda N})$. Clearly, $C(\lambda,N) \sim 1$ for $\lambda N \geq 1$. With $w$ defined as above, we find
\begin{equation*}
\langle f , w(-i \partial_x , N) f \rangle = \frac{1}{2} \big[ \frac{\alpha_1(N,f)}{C(\lambda,N)} - \frac{\alpha_1(N/2,f)}{C(\lambda,N/2)} \big].
\end{equation*}
Suppose that $q\in C^\infty(\mathbb{R} \times \mathbb{T})$ is a solution to \eqref{eq:dNLS}. Let $q_\lambda$ denote the rescaled solution to \eqref{eq:dNLS}:
\begin{equation*}
q_\lambda: \mathbb{R} \times \mathbb{T}_\lambda \to \mathbb{C}, \quad q_\lambda(t,x) = \lambda^{-1/2} q(\lambda^{-2} t, \lambda^{-1} x).
\end{equation*}
With the conventions introduced above, the identities from Lemma \ref{lem:TraceEstimates} and Corollary \ref{lem:LeadingTermCircle} allow for the same error estimates as in the real line case uniformly for $\lambda \in 2^{\mathbb{N}_0}$. We arrive at
\begin{equation*}
\Vert q_\lambda(t) \Vert_{Z^s_r} \lesssim_{r,s} \Vert q_\lambda(0) \Vert_{Z^s_r} + [ \Vert q_\lambda(t) \Vert^{2}_{B_r^{s,2}} + \Vert q_\lambda(0) \Vert^{2}_{B_r^{s,2}} ].
\end{equation*}
By $L^2$-conservation and estimating the $B_{r,2}^{s}$-norm in terms of the $Z^s_r$-norm:
\begin{equation*}
\Vert q_\lambda(t) \Vert_{Z^s_r} \lesssim_{r,s} \Vert q_\lambda(0) \Vert_{Z^s_r} + \Vert q_\lambda(0) \Vert_{L^2_\lambda}^2 + [ \Vert q_\lambda(t) \Vert^{2}_{Z^s_r} + \Vert q_\lambda(0) \Vert^{2}_{Z^s_r} ].
\end{equation*}
As in the real line case, smallness of the $Z^s_r$-norm can be achieved by taking $\lambda \to \infty$ provided that the $L^2$-norm of $q(0)$ is chosen small enough. Also, $\| q_\lambda(0) \|_{L^2(\mathbb{T}_\lambda)} = \| q(0) \|_{L^2}$. Hence, the continuity argument given in the line case proves global a priori estimates in the circle case for small $L^2$-norm of the initial data. The proof of Theorem \ref{thm:APrioriEstimates} is complete.

\section*{Acknowledgements}
Funded by the Deutsche Forschungsgemeinschaft (DFG, German Research Foundation) - Project-ID 258734477 - SFB 1173.
We would like to thank the referee for a very careful reading of an earlier version of the manuscript and insightful comments, which led to some corrections in Sections \ref{section:PerturbationDeterminant} and \ref{section:ConservationSobolevNorms} and helped to improve the presentation.

\end{document}